\documentclass{amsart}
\usepackage{amssymb}
\newtheorem{thm}{Theorem}

\begin{document}

\bibliographystyle{plain}

\title[Milo\v s Arsenovi\'c and Romi F. Shamoyan]{A note on Fefferman-Stein type characterizations for certain spaces of analytic functions on the unit disc}

\author[]{Milo\v s Arsenovi\' c$\dagger$}
\author[]{Romi F. Shamoyan}

\address{Faculty of mathematics, University of Belgrade, Studentski Trg 16, 11000 Belgrade, Serbia}
\email{\rm arsenovic@matf.bg.ac.rs}

\address{Bryansk University, Bryansk Russia}
\email{\rm rshamoyan@yahoo.com}

\thanks{$\dagger$ The first author supported by the Ministry of Science, Serbia, project OI174017}

\date{}

\maketitle

\begin{abstract}

We obtain new characterizations of Bergman and Bloch spaces on the unit disc involving equivalent (quasi)-norms on
these spaces. Our results are in the spirit of estimates obtained by Fefferman and Stein for Hardy spaces in
$\mathbb R^n$.

\end{abstract}

\footnotetext[1]{Mathematics Subject Classification 2010 Primary 30H20, Secondary 30H25, 30H30.  Key words
and Phrases: Hardy spaces, Bergman spaces, equivalent quasinorms, analytic functions.}

\section{Introduction}

We denote by $H(\Omega)$ the space of all analytic functions in a domain $\Omega \subset \mathbb C$, $H^p = H^p(\mathbb D)$ denotes the classical Hardy space on the unit disc $\mathbb D = \{ z \in \mathbb C : |z| < 1\}$ and $A^p = A^p(\mathbb D)$ denotes the Bergman space on $\mathbb D$, see \cite{DS} and \cite{Du}. We set $A^p_0 = \{ f \in A^p : f(0) = 0\}$. Also, $\Gamma_\alpha(\xi)$ denotes the Stolz region at $\xi \in \mathbb T = \{ \xi \in \mathbb C :
|\xi | = 1 \}$ of aperture $\alpha > 1$. For $t > 0$ and an analytic function $f(z) = \sum_{n=0}^\infty a_n z^n$ in the unit disc the fractional derivative of $f$ of order $t$ is defined by $D^t f(z) = \sum_{n=0}^\infty (n+1)^t a_n z^n$, it is also analytic in $\mathbb D$. Area measure on $\mathbb D$ is
denoted by $dm$.
The Bloch space $\mathcal B$, defined by
$${\mathcal B} = \left\{ f \in H(\mathbb D) : \| f \|_{\mathcal B} = \sup_{z \in \mathbb D} |f'(z)|(1-|z|) < \infty \right\}$$
is closely related to Carleson measures and corresponds to the endpoint case $p=1$ of $Q_p$ classes ($0<p\leq 1$),
see \cite{JX1}. Related spaces $\mathcal B_p$, $1 < p < \infty$, are defined by
$${\mathcal B}_p (\mathbb D) = \left\{ f\in H(\mathbb D) : \int_{\mathbb D} |f'(z)|^p (1-|z|)^{p-2} dm(z) = \| f
\|^p_{{\mathcal B}_p} < \infty \right\},$$
note that $\| f \|_{\mathcal B}$ and $\| f \|_{\mathcal B_p}$ are not true norms, but $|f(0)| + \| f \|_{\mathcal B}$
and $|f(0)|+ \| f \|_{\mathcal B_P}$ are norms which make respective spaces into Banach spaces.

The following result is proved by E. G. Kwon in \cite{Kw1}:
\begin{thm}[see \cite{Kw1}]
If $0 < p < \infty$ and $0 \leq \beta < p+2$, then $f \in H(\mathbb D)$ belongs
to $\mathcal B$ if and only if it satisfies the following condition:
\begin{equation}\label{eqA}
\sup_{a \in \mathbb D} \int_{\mathbb D} \int_{\mathbb D} \frac{|f(z) - f(w)|^{p-\beta}}{|1 - \overline w z|^4}
|f'(z)|^\beta (1-|z|)^\beta \frac{(1-|w|)^2(1-|a|)^2}{|1-\overline w a|^4} dm(z) dm(w) < \infty,
\end{equation}
in fact the above expression is equivalent to $\| f \|^p_{\mathcal B}$.

Similarly, if $1 < p < \infty$ and $0 \leq \beta < p+2$, then a function $f \in H(\mathbb D)$ belongs to
${\mathcal B}_p$ if and only if \begin{equation}\label{eqK}
\int_{\mathbb D} \int_{\mathbb D} \frac{|f(w) - f(z)|^{p-\beta}}{|1 - \overline w z|^4} |f'(z)|^\beta
(1 - |z|)^\beta dm(z) dm(w) < \infty,
\end{equation}
moreover, the above expression is equivalent to $\| f \|_{\mathcal B_p}^p$.
\end{thm}

We relate these estimates to the so called Fefferman-Stein type characterizations. By Fefferman-Stein
characterizations we mean the following relations:
\begin{equation}\label{fsr}
\| F \|_X \asymp \inf_{ \omega \in S} \|\Phi( F, \omega) \|_Y,
\end{equation}
where $X$ and $Y$ are (quasi)-normed subspaces of $H(\mathbb D)$, $S = S_F$ is a certain class of measurable functions and $\Phi$ is a nonanalytic function of two variables. This idea was used to determine the predual of $Q_p$ classes,
$0<p\leq 1$, see \cite{JX2}, \cite{JX3}. In certain cases the infimum in (\ref{fsr}) is attained, see \cite{DX}, especially section 5. Using ideas from \cite{FS} and \cite{CMS} the authors there extended and used such characterizations for certain Hardy classes in $\mathbb R^n$. Later one of the authors provided a similar Fefferman-Stein type characterization for the analytic Hardy spaces, this is the second part of the following theorem, the first part is a classical result of N. Lusin.

\begin{thm}[see \cite{Lu},\cite{Sh}]\label{thma}
Let $0 < p, t < \infty$. Then, for $f \in H(\mathbb D)$ we have
\begin{equation}\label{thma1}
\| f \|_{H^p}^p \asymp \int_{\mathbb T} \left( \int_{\Gamma_\alpha(\xi)} |D^t f |^2 (1-|z|)^{2t-2} dm(z)
\right)^{p/2} d\xi.
\end{equation}
Next, let $s > 0$ and $0 < p < 2$. Then, for $f \in H(\mathbb D)$, we have
\begin{equation}\label{thma2}
\int_{\mathbb T} \left( \int_{\Gamma_\alpha(\xi)} |f'(z)|^2 dm(z) \right)^{p/2} d\xi \asymp \inf_{\omega \in S_1}
\left( \int_{\mathbb D} |f'(z)|^s (1-|z|)^{s-1} \frac{dm(z)}{\omega(z)} \right)^{p/2},
\end{equation}
where
$$S_1 = \left\{ \omega \geq 0 : \| \sup_{\Gamma_\alpha(\xi)} \omega(z) (1-|z|)^{2-s} |f'(z)|^{2-s} \|_{L^{
\frac{p}{2-p}}(\mathbb T)} < 1 \right\}.$$
\end{thm}
Here we show that similar results are true for Bloch and Bergman spaces $A^p_0$ in the unit desc. An interesting problem would be to obtain similar results for $Q_p$ or other BMOA-type spaces in the unit disc.

\section{Main results}

In this section we state and prove the main results of this paper. They are analogous to the previously obtained
results on Fefferman-Stein characterizations of Hardy spaces in $\mathbb R^n$ and our proofs heavily rely on the
mentioned results of E. G. Kwon.

\begin{thm}
Let $1 < \alpha < 2$, $1/\alpha + 1/\alpha' = 1$ and $p \geq \alpha'$. Then for $F \in H(\mathbb D)$ with $F(0) = 0$ we have
$$ \| F \|_{A^p_0}^p \asymp \inf_{\omega \in S_2} \left( \int_{\mathbb D} |F'(z)|^\alpha \omega^\alpha (z) dm(z)
\right)^{1/\alpha},$$
where
$$S_2 = \left\{ \omega \geq 0 : \int_{\mathbb D} |F'(z)|^{(p-1)\alpha'} \omega^{-\alpha'}(z)(1-|z|)^{p\alpha'} dm(z) \leq 1 \right\}.$$
\end{thm}

{\it Proof.} Here we use the following result from \cite{JX1}, Chapter 2: If $F \in H({\mathbb D})$ and $F(0) = 0$, then
\begin{equation}\label{ast}
\int_{\mathbb D} |F(z)|^p dm(z) \asymp  \int_{\mathbb D} |F(z)|^{p-\beta} |F'(z)|^\beta (1-|z|)^\beta dm(z),
\end{equation}
where $0<p<\infty$, $0 \leq \beta < p+2$. Taking $ \beta = p$ and applying H\"older inequality we obtain
\begin{align*}
\| F \|_{A^p_0}^p & \leq  \left( \int_{\mathbb D} |F'(z)|^\alpha \omega^\alpha(z) dm(z) \right)^{1/\alpha} \cr
& \times  \left( \int_{\mathbb D} |F'(z)|^{(p-1)\alpha'} \omega^{-\alpha'}(z) (1-|z|)^{p\alpha'} dm(z) \right)^{1/\alpha'}\cr
\end{align*}
and this gives one estimate stated in Theorem 3. It is of some interest to note that this estimate is true
for all $1 < \alpha < \infty$. Now we prove the reverse estimate by choosing a special admissible test function in $S_2$. We can assume $\| F \|_{A^p_0} = 1$ and set
$$\tilde \omega(z) = \frac{|F'(z)|^{p/\alpha}(1-|z|)^{1+p/\alpha}}{|F(z)|}, \qquad z \in \mathbb D.$$
A straightforward calculation shows that
$$\int_{\mathbb D} |F'(z)|^{(p-1)\alpha'} (1-|z|)^{p\alpha'} \tilde\omega^{-\alpha'}(z) dm(z) =
 \int_{\mathbb D} |F'(z)|^{p-\alpha'} (1-|z|)^{p-\alpha'} |F(z)|^{\alpha'} dm(z).$$
Now (\ref{ast}), with $\beta = p -\alpha'$, tells us that the last expression is comparable to $\| F \|_{A^p_0}^p$
and therefore bounded by $C = C_{p, \alpha} > 0$. Hence $\omega = C^{1/\alpha'} \tilde \omega\in S_2$ and we have
\begin{align*}
\int_{\mathbb D} |F'(z)|^\alpha \omega^\alpha(z) dm(z) & =  C \int_{\mathbb D} |F'(z)|^\alpha
\tilde \omega^\alpha(z) dm(z) \cr
& = C \int_{\mathbb D} |F'(z)|^{p+\alpha} |F(z)|^{-\alpha} (1-|z|)^{p+\alpha} dm(z).\cr
\end{align*}
The last integral is, by (\ref{ast}) with $\beta = p + \alpha < p + 2$, bounded from above by a constant depending only on $p$ and $\alpha$ and this ends the proof of Theorem 3.

To simplify formulation of our next theorem we introduce, for $a \in \mathbb D$,
$$d\mu_a(z,w) = \frac{(1-|w|)^2(1-|a|)^2 dm(z)dm(w)}{|1-\overline z w|^4 |1- \overline w a|^4}.$$

\begin{thm}
Let $1 < \alpha < 2$, $1/\alpha + 1/\alpha' = 1$ and $p \geq \alpha'$. Then we have, for $F \in H(\mathbb D)$ with
$F(0) = 0$
\begin{equation}\label{thm2u}
\| F \|_{\mathcal B} \asymp \inf_{\omega \in S_3} \sup_{a \in \mathbb D} \left( \int_{\mathbb D} \int_{\mathbb D}
\omega^\alpha (z, w) |F'(z)|^\alpha d \mu_a (z, w) \right)^{1/\alpha},
\end{equation}
where
\begin{equation}\label{thm2u1}
S_3 = \left\{ \omega \geq 0 : \sup_{a \in \mathbb D} \int_{\mathbb D} \int_{\mathbb D} |F'(z)|^{\alpha'(p-1)}
\omega^{-\alpha'}(z, w) (1-|z|)^{p\alpha'} d \mu_a(z,w) \leq 1 \right\}.
\end{equation}
\end{thm}

{\it Proof.} Let $F \in {\mathcal B}$, setting $\beta = p$ in (\ref{eqA}) we obtain for arbitrary
$\omega \in S_3$:
\begin{align*}
\| F \|_{\mathcal B} & \leq C \sup_{a \in \mathbb D} \int_{\mathbb D} \int_{\mathbb D}
\frac{|F'(z)|^p(1-|z|)^p}{|1-\overline w z|^4} \frac{(1-|w|)^2(1-|a|)^2}{|1-\overline w a|^4} dm(z) dm(w) \cr
& = C \sup_{a \in \mathbb D} \int_{\mathbb D} \int_{\mathbb D} |F'(z)| \omega(z,w) |F'(z)|^{p-1} \omega^{-1}(z,w) (1-|z|)^p d\mu_a(z,w)\cr
& \leq C \sup_{a \in \mathbb D} \left( \int_{\mathbb D}\int_{\mathbb D} |F'(z)|^\alpha \omega^\alpha(z,w) d\mu_a(z,w)
\right)^{1/\alpha}\cr
& \times \sup_{a\in\mathbb D} \left( \int_{\mathbb D}\int_{\mathbb D} |F'(z)|^{\alpha'(p-1)} \omega^{-\alpha'}(z,w)
(1-|z|)^{p\alpha'} d\mu_a(z,w) \right)^{1/\alpha'} \cr
& \leq C \sup_{a \in \mathbb D} \left( \int_{\mathbb D}\int_{\mathbb D} |F'(z)|^\alpha \omega^\alpha(z,w) d\mu_a(z,w)
\right)^{1/\alpha}.\cr
\end{align*}
Taking infimum over all $\omega \in S_3$ one obtains an estimate of $\| F \|_{\mathcal B}$ in terms of the righthand side in (\ref{thm2u}). Now we turn to the reverse estimate, we can assume $\| F \|_{\mathcal B} = 1$. Taking
\begin{equation}\label{test}
\tilde \omega(z, w) = \frac{|F'(z)|^{p/\alpha}(1-|z|)^{\frac{p+\alpha}{\alpha}}}{|f(z) - f(w)|}
\end{equation}
one obtains by an easy calculation and relation (\ref{eqA}) with $\beta = p - \alpha'$
\begin{align*}
& \sup_{a \in \mathbb D} \int_{\mathbb D} \int_{\mathbb D} |F'(z)|^{\alpha' (p-1)} \omega^{-\alpha'}(z, w)
(1-|z|)^{p\alpha'} d\mu_a(z, w)\cr
= & \sup_{a \in \mathbb D} \int_{\mathbb D} \int_{\mathbb D} |F'(z)|^{p-\alpha'}
|f(z) - f(w)|^{\alpha'} (1-|z|)^{p-\alpha'} d\mu_a(z,w) \asymp  \| F \|_{\mathcal B}^p.
\end{align*}
As in the proof of the previous theorem this means that $\omega = C \tilde \omega$ is in $S_3$, where $C = C_{p, \alpha} > 0$. With this choice of $\omega$ we have
\begin{align*}
& \sup_{a \in \mathbb D} \int_{\mathbb D} \int_{\mathbb D} |F'(z)|^\alpha \omega^\alpha (z, w) d\mu_a(z, w) \cr
= & C^\alpha \sup_{a \in \mathbb D} \int_{\mathbb D} \int_{\mathbb D}
|F'(z)|^{p+\alpha} (1-|z|)^{p+\alpha} |f(z) - f(w)|^{-\alpha} d\mu_a (z, w) \cr
\asymp & \| F \|_{\mathcal B}^p = 1,
\end{align*}
where we used (\ref{eqA}) with $\beta = p + \alpha$. This ends the proof of Theorem 4.

This theorem has a counterpart for $\mathcal B_p$ spaces. It is convenient to introduce a measure
$d\mu(z, w) = |1-\overline w z|^{-4} dm(z) dm(w)$.

\begin{thm}
Let $1 < p < \infty$, $1 < \alpha < 2$, $1/\alpha + 1/\alpha' = 1$ and $ p \geq \alpha'$. Then we have, for $F \in
H(\mathbb D)$ with $F(0) = 0$:
\begin{equation}\label{bp}
\| F \|_{\mathcal B_p} \asymp \inf_{\omega \in S_4} \left( \int_{\mathbb D} \int_{\mathbb D} \omega^\alpha (z,w)
|F'(z)|^\alpha d\mu(z,w) \right)^{1/\alpha},
\end{equation}
where
$$S_4 = \left\{ \omega \geq 0 : \int_{\mathbb D} \int_{\mathbb D} |F'(z)|^{\alpha' (p-1)} \omega^{-\alpha'} (z, w)
(1-|z|)^{p\alpha'} d\mu(z,w) \leq 1 \right\}.$$
\end{thm}

{\it Proof.} The proof of this theorem parallels the proof of the previous one. Namely we use (\ref{eqK}) with $\beta = p$ to obtain, for arbitrary $\omega \in S_4$,
\begin{align*}
\| F \|_{\mathcal B_p}^p & \asymp \int_{\mathbb D} \int_{\mathbb D} |F'(z)|\omega(z, w) |F'(z)|^{p-1} (1-|z|)^p
\omega^{-1}(z, w) d\mu(z, w) \cr
& \leq \left( \int_{\mathbb D} \int_{\mathbb D} |F'(z)|^\alpha \omega^\alpha (z, w) d\mu(z,w) \right)^{1/\alpha} \cr
& \times \left( \int_{\mathbb D} \int_{\mathbb D} |F'(z)|^{\alpha' (p-1)} \omega^{-\alpha'}(z, w)
(1-|z|)^{\alpha' p} d\mu(z,w) \right)^{1/\alpha'} \cr
& \leq \left( \int_{\mathbb D} \int_{\mathbb D} |F'(z)|^\alpha \omega^\alpha (z, w) d\mu(z,w) \right)^{1/\alpha}. \cr
\end{align*}
In proving the reverse estimate one can use the same test function as in (\ref{test}), the role of condition
(\ref{eqA}) is taken by condition (\ref{eqK}). We leave details to the reader.

\bibliography{BiblFeffStCh}

\end{document}